\documentclass[conference]{IEEEtran}
\usepackage{times}

\usepackage[numbers]{natbib}
\usepackage{multicol}
\usepackage[bookmarks=true, hidelinks]{hyperref}
\usepackage{amsmath}
\usepackage[capitalize]{cleveref}
\usepackage{graphicx}
\usepackage{stfloats}
\usepackage[acronym]{glossaries}
\newcommand{\myparagraph}[1]{\noindent\textbf{#1}}
\newacronym{acr:amod}{AMoD}{Autonomous Mobility-on-Demand}
\newacronym{acr:av}{AV}{autonomous vehicles}
\newacronym{acr:milp}{MILP}{Mixed-Integer Linear Program}
\newacronym{acr:ilp}{ILP}{Integer Linear Program}
\newacronym{acr:vrp}{VRP}{Vehicle Routing Problem}
\newacronym{abk:av}{AV}{Autonomous Vehicle}
\renewcommand{\baselinestretch}{0.96}

\begin{document}

\title{Robo-Taxi Fleet Coordination with Accelerated High-Capacity Ridepooling}

\author{\authorblockN{Xinling Li\authorrefmark{1}, Daniele Gammelli\authorrefmark{2}, Alex Wallar\authorrefmark{3}, Jinhua Zhao\authorrefmark{4}, Gioele Zardini\authorrefmark{1}}
\authorblockA{\authorrefmark{1}Laboratory for Information \& Decision Systems, 
Massachusetts Institute of Technology}
\authorblockA{\authorrefmark{2}Department of Aeronautics and Astronautics, Stanford University}
\authorblockA{\authorrefmark{3}The Routing Company}
\authorblockA{\authorrefmark{4}Department of Urban Studies and Planning, Massachusetts Institute of Technology}
}


%

\maketitle

\begin{abstract}
Rapid urbanization has led to a surge of customizable mobility demand in urban areas, which makes on-demand services increasingly popular. On-demand services are flexible while reducing the need for private cars, thus mitigating congestion and parking issues in limited urban space. While the coordination of high-capacity ridepooling on-demand service requires effective control to ensure efficiency, the emergence of the paradigm of robo-taxi opens the opportunity for centralized fleet control for an improved service quality. In this work, we propose two acceleration algorithms for the most advanced large-scale high-capacity algorithm proposed in ~\cite{alonso2017demand}. We prove the improvement in the real-time performance of the algorithm by using real-world on-demand data from Manhattan, NYC.
\end{abstract}

\IEEEpeerreviewmaketitle

\section{Introduction} 
By 2050, 68\% of the global population will reside in urban areas~\cite{un_wup_2018}.
Increasing urbanization has intensified demand for customized mobility solutions, contributing to worsening congestion and greenhouse gas emissions.
In response, on-demand mobility services have emerged as a promising alternative to private vehicle ownership, offering flexibility while potentially alleviating environmental and infrastructural stress.
A central challenge in realizing the full potential of such technologies lies in the design of effective vehicle coordination strategies.
Indeed, due to the inherently spatiotemporally unbalanced nature of travel demand, poor coordination can result in significant empty mileage and deadhead trips, ultimately diminishing the intended efficiency and sustainability gains~\cite{oh2021impacts}.

The advent of \gls{acr:amod} systems, also referred to as robo-taxi services, introduces the concept of centralized, algorithmically driven fleet control~\cite{zardini2022analysis}.
At the core of such systems is the assignment of vehicles to service requests, a task that has attracted considerable attention in the \gls{acr:amod} literature.
Prior approaches span optimization-based model predictive control~\cite{carron2019scalable,tsao2018stochastic}, learning-based methods~\cite{sadeghi2022reinforcement,han2016routing}, and hybrid frameworks integrating optimization with learning~\cite{liang2021integrated,jungel2023learning,woywood2024multi,tresca2025robo}.
However, these contributions predominantly focus on single-passenger assignments.
The case of \gls{acr:amod} with ridepooling, despite its potential to further reduce congestion and emissions by lowering fleet size and merging trips, remains underexplored.
\begin{figure}[tb]
    \centering
    \includegraphics[width=\linewidth]{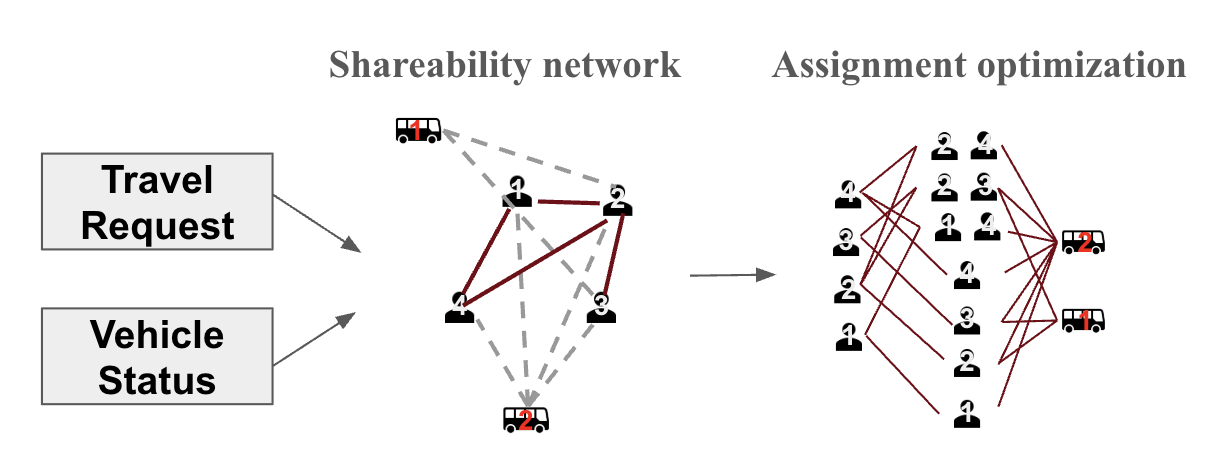}
    \caption{Assignment algorithm for each decision window.}
    \label{fig:original}
\end{figure}
\myparagraph{Related work} --
Within the literature on \gls{acr:amod} coordination with sharing, much of the focus has been on low-capacity ridepooling, typically allowing at most two passengers to share a ride~\cite{enders2023hybrid,jindal2018optimizing,tsao2018stochastic,yu2019integrated,zardini2022analysis}.
While this restriction simplifies the coordination problem, it undermines one of the primary advantages of ridepooling, namely, \emph{sharing}, thus limiting overall system efficiency~\cite{alonso2017demand}.
In contrast, high-capacity ridepooling promises substantial gains in efficiency and sustainability by accommodating more passengers per vehicle.
However, its associated combinatorial complexity makes real-time optimization in large-scale, real-world networks computationally intractable using purely optimization-based approaches.

High-capacity ridepooling research involves a fundamental trade-off between \emph{computational efficiency} and \emph{solution quality}.
To address this, some approaches incorporate learning-based modules to accelerate optimization while preserving feasibility guarantees~\cite{shah2020neural,kim2024learning}, whereas others rely on heuristics that simplify sharing scenarios to reduce complexity~\cite{tong2018unified,simonetto2019real}.
However, such heuristics often rely on restrictive assumptions, limiting their applicability in diverse operational settings.

The leading algorithm for high-capacity ridepooling~\cite{alonso2017demand} improves tractability by using a shareability network~\cite{santi2014quantifying} to prune the search space, thus reducing the problem dimension while preserving optimality guarantees.
Its modular \emph{architecture} has made it a cornerstone for follow-up work~\cite{simonetto2019real,alonso2017predictive,cap2018multi}.
However, it still falls short of real-time performance, highlighting the need for further efficiency improvements.

\myparagraph{Contribution} -- 
We propose two algorithmic accelerations based on~\cite{alonso2017demand}, improving its real-time performance.
We validate our approach using a large-scale real-world ridepooling dataset from Manhattan, NYC, and demonstrate improvements in computational efficiency without compromising performance.
\section{Methodology}
Before introducing our acceleration strategies, we first present a concise overview of the original algorithm, illustrated in \cref{fig:original}.
The system assumes a fleet of \glspl{acr:av} coordinated to serve on-demand ride requests, where each request specifies a pickup and dropoff time window, defined by a maximum allowable waiting time and delay.
When starting each decision window, the algorithm aggregates all active requests and computes a \emph{shareability graph} based on the current set of requests and the status of the vehicle fleet (i.e., locations and available capacities).
Subsequently, the algorithm enumerates feasible trips for each vehicle using the shareability graph and formulates a global assignment problem as a \gls{acr:milp}.
Building the shareability graph is computationally efficient, as it is fully parallelizable over the sets of \glspl{acr:av} and requests.
In contrast, the principal computational bottleneck lies in the formulation of the assignment \gls{acr:ilp}.
Specifically, adding an edge to the assignment graph requires solving a \gls{acr:vrp} to evaluate the feasibility of each candidate trip.
Although individual \glspl{acr:vrp} are relatively small, bounded by vehicle capacity, one \gls{acr:vrp} must be solved for every clique in the shareability graph. As the number of cliques grows exponentially with the graph's size, and the enumeration of cliques is inherently sequential with respect to clique size, this process becomes computationally intensive and non-parallelizable. To address this challenge and improve the algorithm's real-time performance, we propose two acceleration techniques described in the following sections.
%

\myparagraph{Data-driven \gls{acr:ilp} construction} --
The original algorithm operates in a purely model-based manner, without leveraging historical data.
To accelerate the construction of the assignment \gls{acr:ilp}, we introduce a data-driven \emph{feasibility prediction layer} that precedes the evaluation of each \gls{acr:vrp}.
Specifically, we learn to predict the feasibility of a candidate clique using prior data, and only solve the full \gls{acr:vrp} if the predicted probability of feasibility exceeds a given threshold.
This allows us to bypass the computational cost of solving \glspl{acr:vrp} for low-likelihood cliques, focusing resources on more promising ones.
The feasibility prediction model is implemented leveraging a Structure2Vec graph embedding~\cite{khalil2017learning}, which encodes the structural properties of each clique to inform the prediction.

\myparagraph{\textbf{Shareability graph partition}} --
The primary source of computational complexity in the assignment \gls{acr:ilp} lies in the density and size of the shareability graph. 
To mitigate this, we propose a \emph{partitioning strategy} whereby the shareability graph is \emph{decomposed} into smaller, disjoint subgraphs.
Feasible trip generation is then performed in parallel within each subgraph.
The resulting subproblems are subsequently aggregated into a unified \gls{acr:ilp}, which follows the original formulation.
This approach not only unlocks parallelization across subgraphs, but also reduces the computational complexity of each instance, thereby significantly improving overall scalability.
%
\section{Numerical case study}
We evaluate the proposed acceleration methods against the original algorithm (baseline) using the NYC TLC dataset~\cite{nyc_tlc_trip_data}, which captures real-world on-demand mobility.
Experiments focus on the morning peak (7:00-9:00 AM, May 15, 2024), comprising ca. 26,000 trips on a network of roughly 4,000 nodes and 10,000 edges.
Decisions are updated every minute.

In real-time deployments, solving the assignment problem to optimality within each decision window is often impractical.
A common approach is to enforce a fixed time budget and deploy the best solution found within that limit.
To reflect this, we apply a 30-second timeout in all experiments, reserving time for communication and dispatch within the 1-minute window.
This setup ensures realism and fairness in comparing the baseline and accelerated methods.

\begin{figure}[tb]
    \centering
    \includegraphics[width=0.95\linewidth]{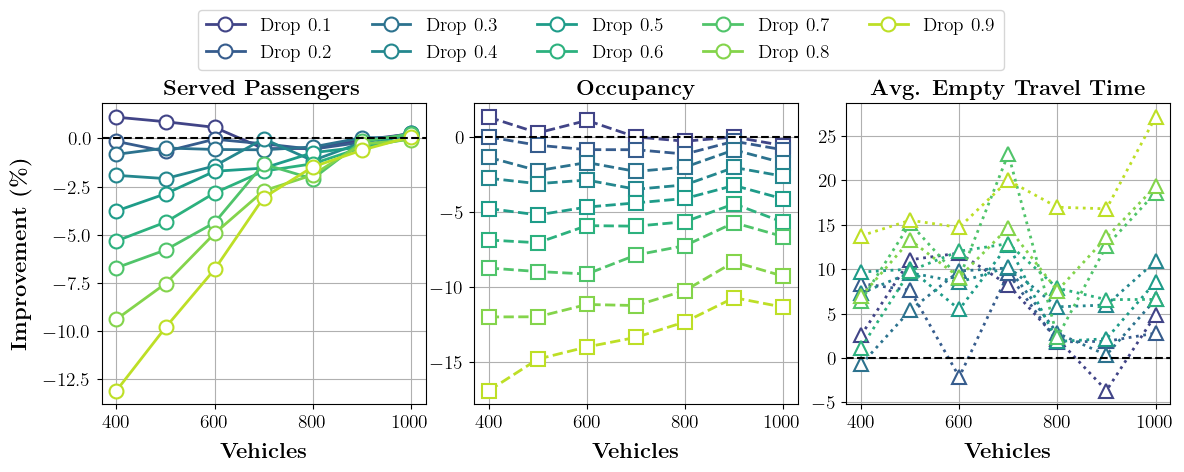}
    \caption{Comparison under different drop rates.}
    \label{fig:data-driven}
\end{figure}

\myparagraph{Data-driven ILP construction} --
The data-driven method selectively skips cliques to prioritize evaluation of larger, more impactful ones.
A common baseline heuristic is to randomly drop cliques during shareability graph construction.
As shown in \cref{fig:data-driven}, the data-driven approach matches the best-performing drop rate in terms of passengers served and average occupancy, while reducing empty travel time.


\begin{figure}[t]
    \centering
    \includegraphics[width=0.95\linewidth]{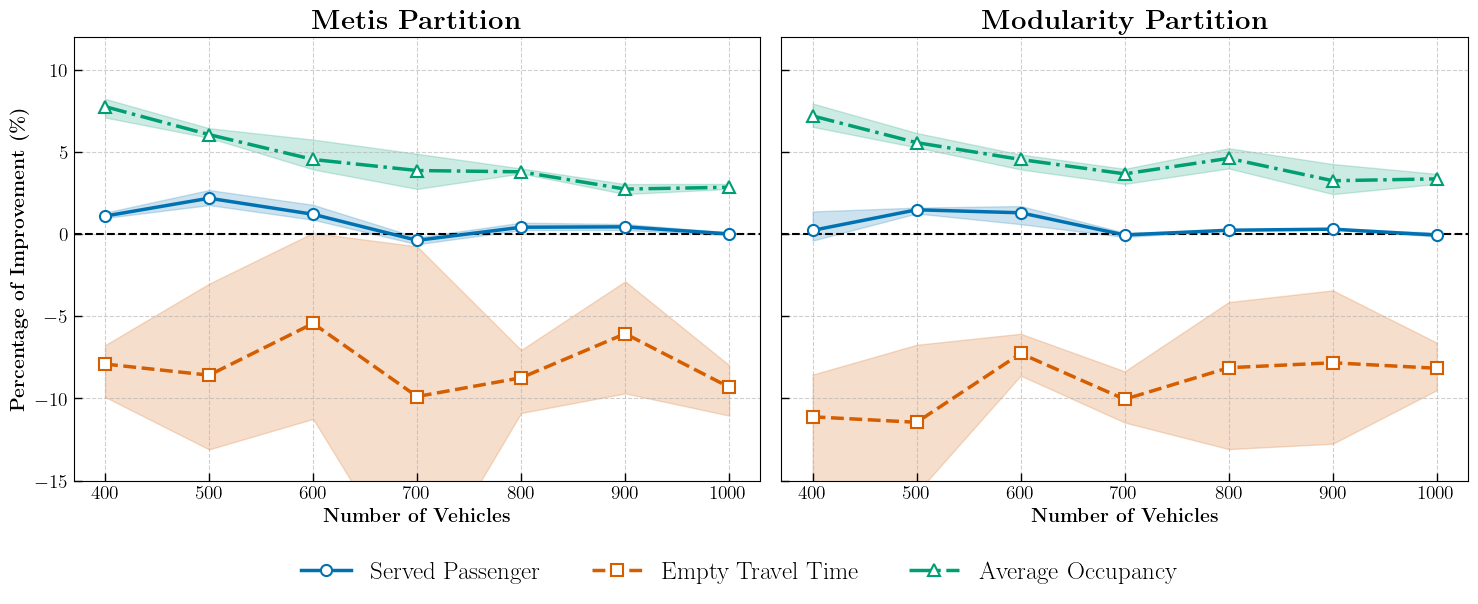}
    \caption{Performance comparison with partition.}
    \label{fig:partition}
\end{figure}

\myparagraph{Shareability graph partition} --
We partition the shareability graph using the multilevel~$k$-way algorithm~\cite{karypis1998multilevelk}, which balances subgraph sizes while minimizing inter-partition edges to enable efficient parallelization.
For comparison, we also test modularity-based partitioning~\cite{clauset2004finding}, which favors edge density over balance.
We can assess the trade-offs between partition quality and load balancing for parallel trip generation. 
\Cref{fig:partition} shows that with partitioning, a higher occupancy and lower empty traveling time is achieved.

\section{Conclusion and future work}
Experimental results show that both acceleration methods substantially improve baseline performance under realistic time constraints, by reducing empty travel time, a key indicator of performance.
This has practical benefits for sustainability and congestion in dense urban settings like NYC.
As future work, we will evaluate accelerations in various demand scenarios to understand their effectiveness under different conditions, and explore the co-design of fleet sizes and composition~\cite{zardini2022co}.



\section{Acknowledgement}
This work is supported by the US Department of Energy, Grant Agreement DE-EE0011186.

\bibliographystyle{plainnat}
\bibliography{references}

\end{document}